\documentclass[
    ,final            
  ]
{aipproc}
\usepackage{amssymb, bm}
\usepackage[margin=2.5cm]{geometry}
\layoutstyle{8x11single}

\newtheorem{defi}{Definition}

\newtheorem{rem}{Remark}
\newtheorem{theo}{Theorem}

\newtheorem{lem}{Lemma}

\def\RR{\mathbb{R}}

\def\pmatrix{\left(\begin{array}}
\def\endpmatrix{\end{array}\right)}
\def\dd{{\mathrm d}}

\def\bfb{{\bf b}}
\def\bfc{{\bf c}}

\def\aa{{\alpha}}

\def\balfa{\bm{\alpha}}
\def\tgam{\bar{\gamma}}
\def\hgamma{\hat{\gamma}}

\def\II{{\cal I}}
\def\PP{{\cal P}}
\def\dd{{\mathrm{d}}}

\def\diag{{\rm diag}}

\def\heta{{\hat\eta}}

\def\tgam{{\tilde{\gamma}}}

\begin{document}

\title{Enhanced HBVMs for the numerical solution of Hamiltonian problems with multiple invariants}


\author{Luigi Brugnano}{
  address={Dipartimento di Matematica e Informatica ``U.\,Dini'', Universit\`a di Firenze, Italy}
}

\author{Yajuan Sun}{
  address={Academy of Mathematics and Systems Science, Chinese Academy of Sciences, Beijing, China}
}

\begin{abstract} Recently, the class of energy-conserving Runge-Kutta methods named Hamiltonian Boundary Value Methods (HBVMs), has been proposed for the efficient solution of Hamiltonian problems, as well as for other types of conservative problems. In this paper, we report further advances concerning such methods, resulting in their  {\em enhanced} version  (Enhanced HBVMs, or EHBVMs). The basic theoretical results are sketched, along with a few numerical tests on a Hamiltonian problem, taken from the literature, possessing multiple invariants.
\end{abstract}

\maketitle


\section{INTRODUCTION}

In a recent series of papers \cite{BIT09,BIT09_1,BIS10,BIT10,BIT11,BIT12,BIT12_1} (see also \cite{BI12_1,BI12_2}), the class of Hamiltonian Boundary Value Methods (HBVMs) has been proposed for the efficient numerical solution of Hamiltonian problems, i.e., problems in the form
\begin{equation}\label{ham}
y' = J\nabla H(y),  \qquad  y\in \RR^{2m},
\end{equation}
where $y=(q,p)^T$, $J = \pmatrix{cc} 0 & I_m\\ -I_m &0\endpmatrix = -J^T = -J^{-1}$, and $H(y)\equiv H(q,p)$ is the {\em Hamiltonian} (also referred to as the {\em energy}).  The system turns out to be {\em energy-preserving}. Indeed, since $J$ is   skew-symmetric, from (\ref{ham}) one has that
$$\frac{\dd}{\dd t} H(y(t)) = \nabla H(y(t))^Ty'(t) = \nabla H(y(t))^TJ\nabla H(y(t)) = 0,$$
which implies that $H(y(t))\equiv H(y(0))$ for all $t\ge0$. The basic idea on which HBVMs rely is that of {\em discrete line integral}, introduced in \cite{IP07,IP08,IT09}. Such a tool has also been extended to different kinds of conservative problems \cite{BCMR12,BI12} (see also \cite{TS12}). We here use the approach of discrete line integrals, to derive the {\em enhanced} version of HBVMs, able to cope with problems, in the form (\ref{ham}), possessing additional (functionally independent) invariants, besides the Hamiltonian. In more details, let us assume  that
\begin{equation}\label{nuinv}
L:\RR^{2m}\rightarrow \RR^\nu
\end{equation}
is a set of smooth invariants of  system (\ref{ham}). Then, one has
$$\nabla L(y)^T J\nabla H(y) = 0\in\RR^\nu, \qquad \forall y\in \RR^{2m},$$
where $\nabla L(y)^T$ is the Jacobian of $L$. The {\em ehanced} version of HBVMs, which we call {\em Enhanced HBVMs (EHBVMs)}, will be able, under suitable mild hypotheses, to provide a discrete solution for which all the invariants (\ref{nuinv}) are conserved, besides the Hamiltonian. Full theoretical details on EHBVMs can be found in \cite{BS13}.

\section{Mutiple invariants preserving HBVMs}

Let $\sigma$ be a polynomial of  degree $s$, approximating the solution of (\ref{ham}) on the interval $[0,h]$,  in the form
\begin{equation}\label{sigma1}
\sigma'(ch) = \sum_{j=0}^{s-1} P_j(c) \gamma_j(\sigma), \qquad c\in[0,1],\qquad \mbox{with}\quad  s>\nu,
\end{equation}
where $\{P_j\}_{j\ge0}$ is the family of normalized and shifted Legendre polynomials,  orthonormal on the interval $[0,1]$,
i.e., $$P_j\in\Pi_j,\qquad \int_0^1 P_i(x)P_j(x)\dd x = \delta_{ij}. \qquad \forall i,j\ge0.$$
We require that vector coefficients $\gamma_j(\sigma)$ in (\ref{sigma1})  satisfy
\begin{equation}\label{gammaj}
\gamma_j(\sigma) = \eta_jJ\int_0^1P_j(\tau) \nabla H(\sigma(\tau h))\dd\tau\equiv \eta_j\tgam_j(\sigma), \qquad \mbox{with}\qquad \eta_j\in\RR,  \qquad  j=0,\dots,s-1.
\end{equation}
Imposing $\sigma(0)=y_0$, and defining the new approximation as $y_1\equiv\sigma(h)\approx y(h)$, provides energy-conservation, since due to the skew-symmetry of matrix $J$, one obtains \cite{BIT10,BIT12_1}:
\begin{eqnarray*}
H(y_1)-H(y_0) &=& H(\sigma(h))-H(\sigma(0)) ~=~ \int_0^h \nabla H(\sigma(t))^T\sigma'(t)\dd t\\
 &=& h\int_0^1 \nabla H(\sigma(\tau h))^T\sigma'(\tau h)\dd\tau ~=~h\sum_{j=0}^{s-1}\eta_j \tgam(\sigma)_j^TJ\tgam_j(\sigma)~=~0. 
 \end{eqnarray*}
Taking $\eta_i=1$, for $0\leq i\leq s-1$, provides the energy-conserving methods named  HBVMs \cite{BIT10}. Here we set, instead,
\begin{equation}\label{coe}
\eta_j=1, \quad 0\leq j\leq s-\nu-1,\qquad \eta_j = 1-h^{2(s-1-j)}\aa_j,\quad s-\nu\leq j\leq s-1,
\end{equation} with the coefficients $\aa_j$ determined by imposing the conservation of the additional invariants (\ref{nuinv}) at $y_1$:
\begin{eqnarray*}
 L(y_1)-L(y_0) &=& L(\sigma(h))-L(\sigma(0))~=~h\sum_{j=0}^{s-1} \left[\int_0^1P_j(\tau)\nabla L(\sigma(\tau h))\dd\tau\right]^T\gamma_j(\sigma)\dd\tau \\
&\equiv& h\left[\sum_{j=0}^{s-1}\phi_j(\sigma)^T\tgam_j(\sigma) ~-~\sum_{j=s-\nu}^{s-1} h^{2(s-1-j)}\aa_j\phi_j(\sigma)^T\tgam_j(\sigma)\right],
\end{eqnarray*}
where
\begin{equation}\label{fij}
\phi_j(\sigma) = \int_0^1P_j(\tau)\nabla L(\sigma(\tau h))\dd\tau ~\in\RR^{2m\times\nu}, \qquad j\ge0,
\end{equation}
and $\{\tgam_j(\sigma)\}_{j=0}^{s-1}$ is defined according to (\ref{gammaj}).
Consequently, ~$L(y_1)=L(y_0)$~ if and only if
\begin{equation}\label{Lcons}
\sum_{j=s-\nu}^{s-1} h^{2(s-1-j)}\aa_j\phi_j(\sigma)^T\tgam_j(\sigma) = \sum_{j=0}^{s-1}\phi_j(\sigma)^T\tgam_j(\sigma).
\end{equation}
By defining the matrix
$$\Gamma(\sigma) = \left[\begin{array}{ccc} h^{2(\nu-1)}\phi_{s-\nu}(\sigma)^T\tgam_{s-\nu}(\sigma), & \dots ~,& h^0\phi_{s-1}(\sigma)^T\tgam_{s-1}(\sigma)\end{array}\right]\in\RR^{\nu\times\nu},$$ and the vectors
$$\balfa = \left[\begin{array}{ccc} \aa_{s-\nu}, & \dots~, & \aa_{s-1}\end{array}\right]^T,\quad \bfb(\sigma) = \sum_{j=0}^{s-1}\phi_j(\sigma)^T\tgam_j(\sigma)\quad\in\quad\RR^\nu,$$
 equation (\ref{Lcons})  can be recast in vector form as
\begin{equation}\label{Lcons1}
\Gamma(\sigma) \balfa = \bfb(\sigma).
\end{equation}

The following results hold true \cite{BS13}.

\begin{lem}\label{bfb}
 With reference to (\ref{Lcons1}), one has that  ~$\bfb(\sigma) = O(h^{2s})$~ and  ~$\Gamma(\sigma)=O(h^{2s-2})$.
\end{lem}

\begin{theo}\label{alfa2} Assume that matrix $\Gamma(\sigma)$ is nonsingular. Then, the vector $\balfa$ in (\ref{Lcons1}) has $O(h^2)$ entries. The polynomial approximation $\sigma$ defined by (\ref{sigma1}) with
$\eta_j$ in form of (\ref{coe})  conserves all the invariants and, moreover,  $$\sigma(h)-y(h)=O(h^{2s+1}).$$\end{theo}

\section{Discretization and EHBVM$(k,s)$ methods}
 It is noticed that the previous formulae provide an effective method only after that the integrals in (\ref{gammaj}) and (\ref{fij}) are approximated by using a suitable quadrature formula.  For this purpose, we  choose the numerical integration formula defined at the $k$ Gauss-Legendre points $0<c_1<\cdots<c_k<1$.  This leads to  a new polynomial approximation of degree $s$, say $u$, defined by
\begin{eqnarray}\nonumber
u(ch) &=& y_0 + h\sum_{j=0}^{s-1} \int_0^cP_j(x)\dd x \, \heta_j \hgamma_j\\\label{uu}
 &\equiv& y_0 + h\left[\sum_{j=0}^{s-1} \int_0^cP_j(x)\dd x\, \hgamma_j -  \sum_{j=s-\nu}^{s-1}  \int_0^cP_j(x)\dd x\, h^{2(s-1-j)}\hat\aa_j \hgamma_j\right], \qquad c\in[0,1],
\end{eqnarray}
where  $\hgamma_j$, $\heta_j$, and $\hat\aa_j$ are the approximations to $\gamma_j$, $\eta_j$,  and $\aa_j$, respectively, obtained by using the given quadrature formula $(c_\ell,b_\ell)_{\ell=1}^k$.
Observe that $u(0)=y_0$. Setting $y_1\equiv u(h)$ and $u_\ell=u(c_\ell h)$, $\ell=1,\dots,k$, from (\ref{uu}) one has:
\begin{equation}\label{y}
y_1\equiv u(h) = y_0 + h\hgamma_0 = y_0 + h\sum_{\ell=1}^k b_\ell J\nabla H(u_\ell).
\end{equation}
\begin{defi}
We name Enhanced HBVM($k$, $s$) (in short, EHBVM$(k,s)$), the methods defined by (\ref{uu})-(\ref{y}). \end{defi}

The following results hold true \cite{BS13}.

\begin{theo}\label{ordu1} Assuming that both $H$ and $L$ are suitably regular, for all $k\ge s$ the numerical solution generated by a EHBVM$(k,s)$ method satisfies $$y_1-y(h)=O(h^{2s+1}).$$ That is, the method has order $2s$.\end{theo}

\begin{theo}\label{exactL} Assume that the invariants (\ref{nuinv}) of problem (\ref{ham}) are polynomials of degree less than or equal to $\mu=\lfloor2k/s\rfloor$. Then, an EHBVM$(k,s)$ method is invariants-conserving. Moreover, for all general and suitably regular $L$, one has $$L(y_1)-L(y_0) = O(h^{2k+1}).$$\end{theo}

\begin{rem} As a consequence, even though EHBVM$(k,s)$ has order $2s$, one can recover a {\em practical} invariant-conservation (i.e., to within machine round-off), provided that $k$ is large enough.\end{rem}

We end this section by stating the following Runge-Kutta type formulation of a EHBVM$(k,s)$ method \cite{BS13}:
$$\begin{array}{c|c}
\bfc & \II_s\Sigma_s\PP_s^T\Omega\\
\hline &\bfb^T\end{array},$$
where, as usual, ~$\bfc=(c_1,\dots,c_k)^T$, ~$\bfb=(b_1,\dots,b_k)^T$~ are the abscissae and weights vectors, respectively, and
$$\PP_s = \left( P_{j-1}(c_i)\right),~\II_s = \left( \int_0^{c_i} P_{j-1}(x)\dd x\right)\in\RR^{k\times s},\quad\Sigma_s = \diag(1,\heta_1,\dots,\heta_{s-1}),\quad\Omega = \diag(b_1,\dots,b_k).$$
In particular, when $k=s$ and $\heta_j=1$, $j=1,\dots,s-1$, one retrieves the $s$-stage Gauss method.

\section{Numerical tests} We consider the problem defined by the Hamiltonian (\cite{MP04}, see also \cite{BIT12_3})
\begin{equation}\label{macie}
H(q,p) = \frac{1}2 p^Tp + \left( q^Tq \right)^2, \qquad q,p\in\RR^2,\end{equation} admitting the angular momentum,
~$L(q,p) = q_1p_2-q_2p_1$,~  as a further invariant.
In Table~\ref{testres} we list the obtained numerical results, in terms of conservation of the invariants and solution error, by using the following 4-th order methods: the symplectic 2-stages Gauss method which, evidently, will preserve the angular momentum but not the energy; the energy-conserving HBVM(4,2) method which, however, will not conserve the angular momentum;  the fully-conserving EHBVM(4,2) method, conserving both invariants. The obtained results clearly confirm the effectiveness of the new methods.

\begin{table}[t]
\caption{Numerical result for problem (\ref{macie}), with initial point is $q^0 = (1,~1)^T$, $p^0=(10^{-1},~0)^T$, and integration interval $[0,10^2]$; $e_H$ is the Hamiltonian error; $e_L$ is the error in the angular momentum; $e_{sol}$ is the error in the computed solution.}\label{testres}
\begin{tabular}{|r|ccc|ccc|ccc|}
\hline
       &\multicolumn{3}{c|}{2-stage Gauss} &\multicolumn{3}{c|}{HBVM(4,2)} &\multicolumn{3}{c|}{EHBVM(4,2)}\\
$h$ & $e_H$ & $e_L$  &$e_{sol}$ & $e_H$ & $e_L$ &$e_{sol}$ & $e_H$ & $e_L$ &$e_{sol}$\\
\hline
   $10^{-1}$                  &  2.05e-04 &  6.25e-16 &  1.08e-02 & 4.44e-15  &   8.86e-07 & 7.17e-03  & 5.20e-14 &  1.53e-15 & 2.36e-03 \\
   $2^{-1}\cdot10^{-1}$ & 1.26e-05  &  9.71e-16 &  6.83e-04 & 1.87e-14  &   5.55e-08 &  4.55e-04 & 4.53e-14 &  1.19e-15 & 1.51e-04\\
   $2^{-2}\cdot10^{-1}$ &  7.82e-07 &  1.47e-15 &  4.28e-05 & 7.11e-15  &   3.47e-09 & 2.86e-05 & 4.26e-14 &  1.14e-15 & 9.50e-06 \\
   $2^{-3}\cdot10^{-1}$ &  4.88e-08 &  1.42e-15 &  2.67e-06 & 1.07e-14  &   2.17e-10 & 1.79e-06 & 2.04e-14 &  2.64e-15 & 5.95e-07 \\
   $2^{-4}\cdot10^{-1}$ &  3.05e-09 &  2.75e-15 &  1.67e-07 & 9.77e-15  &   1.36e-11 & 1.12e-07 & 1.42e-14 &  3.64e-15 &  3.72e-08\\
\hline
\end{tabular}
\end{table}

\bibliographystyle{aipproc}

\begin{theacknowledgments} The second author was supported by  the Foundation for Innovative
Research Groups of the NNSFC (11021101).
\end{theacknowledgments}

\end{document}